\documentclass[10pt,bezier]{article}
\usepackage{amsmath,graphicx,amssymb,amsfonts,color}

\font\bg=cmbx10 scaled\magstep1 
\font\small=cmr8

\newtheorem{newlemma}{{\bf Lemma}}

\newtheorem{newteorem}{{\bf Theorem}}

\newenvironment{teorem}{\begin{newteorem}{\hspace{-0.5
em}{\bf.}}}{\end{newteorem}}

\newtheorem{newkorolari}{{\bf Corollary}}

\newenvironment{korolari}{\begin{newkorolari}{\hspace{-0.5
em}{\bf.}}}{\end{newkorolari}}

\newtheorem{newkonjek}{{\bf Conjecture}}

\newtheorem{newdefine}{{\bf Definition}}

\newenvironment{define}{\begin{newdefine}{\hspace{-0.5
em}{\bf.}}}{\end{newdefine}}

\begin{document}
\tolerance=10000 \baselineskip18truept
\newbox\thebox
\global\setbox\thebox=\vbox to 0.2truecm{\hsize
0.15truecm\noindent\hfill}
\def\boxit#1{\vbox{\hrule\hbox{\vrule\kern0pt
     \vbox{\kern0pt#1\kern0pt}\kern0pt\vrule}\hrule}}
\def\qed{\lower0.1cm\hbox{\noindent \boxit{\copy\thebox}}\bigskip}
\def\ss{\smallskip}
\def\ms{\medskip}
\def\nt{\noindent}

 \vspace{.3cm}

\centerline {\large\bf More on domination polynomial and domination root}
\bigskip

\bigskip
\baselineskip12truept \centerline{ Saeid
Alikhani$^{a,b,}$\footnote{\baselineskip12truept\it\small E-mail:
alikhani@yazduni.ac.ir} and Emeric Deutsch$^{c}$ }
\baselineskip20truept \centerline{ $^a$Department of Mathematics,
Yazd University} \vskip-8truept \centerline{ 89195-741, Yazd,
Iran}

\centerline{ $^b$School of Mathematics, Institute for
Research in Fundamental Sciences (IPM)} \vskip-8truept\centerline{P.O. Box:
19395-5746, Tehran,
             Iran.}

 \centerline{
$^{c}$  Polytechnic Institute of New York University, United States}
 \nt\rule{12cm}{0.1mm}

\nt{\bg ABSTRACT}

\medskip
\nt Let $G$ be a simple graph of order $n$.  The domination
polynomial of $G$ is the polynomial $D(G,\lambda)=\sum_{i=0}^{n}
d(G,i) \lambda^{i}$, where $d(G,i)$ is the number of dominating
sets  of $G$ of size $i$. Every root of $D(G,\lambda)$ is called
 the domination root of $G$. It is clear that $(0,\infty)$ is zero free interval for
domination polynomial of a graph. It is interesting  to
investigate graphs which have complex domination roots with
positive real parts.
  In this paper,   we first investigate complexity of the domination polynomial at specific points.  Then we present and investigate some families of
graphs whose complex domination roots have positive real part.

\nt{\bf  Mathematics Subject Classification:} {\small 05C60.}

\nt\rule{14cm}{0.1mm}

\section{Introduction}

\nt Let $G$ be a simple graph. For any vertex $v\in V$, the {\it
open neighborhood} of $v$ is the set $N(v)=\{u \in V|uv\in E\}$
and the {\it closed neighborhood}
 is the set $N[v]=N(v)\cup \{v\}$.
 For a set $S\subseteq V$, the open neighborhood of $S$ is $N(S)=\bigcup_{v\in S} N(v)$ and the
  closed neighborhood of $S$ is $N[S]=N(S)\cup S$.
A set $S\subseteq V$ is a {\it dominating set} if $N[S]=V$, or
equivalently,
 every vertex in $V\backslash S$ is adjacent to at least one vertex in $S$.
  An {\it $i$-subset} of $V(G)$ is a subset of
$V(G)$ of cardinality $i$.
 Let ${\cal D}(G,i)$ be the family of
 dominating sets of $G$ which are $i$-subsets and let $d(G,i)=|{\cal D}(G,i)|$.
 The polynomial $D(G,x)=\displaystyle\sum_{i=0}^{|V(G)|} d(G,i) x^{i}$ is defined as {\it domination polynomial} of $G$ (\cite{euro,ars}).
 A root of $D(G,x)$
is called a {\it domination  root} of $G$. We denote the set of
all roots of $D(G,x)$ by $Z(D(G,x))$. For more information and
motivation of domination polynomial  and domination roots refer to
\cite{contem,euro,ars,brown}.

\nt The value of a graph polynomial at a specific point can give sometimes a surprising information
about the structure of the graph. Balister, et al. in \cite{interlace} proved that for any
graph $G$, $|q(G,-1)|$ is always a power of $2$, where $q(G,x)$ is interlace polynomial of a
graph $G$. Stanley in \cite{stanley} proved that $(-1)^nP(G,-1)$ is the number of acyclic orientations
of $G$, where $P(G,x)$ is the chromatic polynomial of $G$ and $n =|V(G)|$. Alikhani in \cite{gcom}
studied the domination polynomial at $-1$ and gave a construction showing that for each
odd number $n$ there is a connected graph $G$ with $D(G,-1) = n$.
Obviously $(0,\infty)$ is zero free interval for domination polynomial of a graph. It is interesting
to find and study graphs with domination roots in the right-half plane.
In this paper we show that the computation of values of $D(G,\lambda)$ for $\lambda \in \mathbb{Q} \setminus \{-2,-1, 0\}$
is \#P-hard. Also present and study graphs which have complex domination roots with
positive real parts.

\nt
In Section 2 we investigate complexity of the domination polynomial at specific points. In
Section 3 we consider specific families of graphs and compute their domination polynomials.
We show that no nonzero real number is domination root of these kind of graphs.
Also we see that these kind of graphs have domination roots in the right half-plane. In
Section 4 we study the domination polynomial of other classes of graphs.

\section{Complexity of domination polynomial at specific points}

\nt
In this section we investigate Turing complexity of the domination polynomial. First we
recall a formula for computing the domination polynomial of the following graph composition.
Let $G$ and $H$ be graphs, with $V(G) = \{v_1,\ldots, v_n\}$. The graph $G\diamond H$ formed by substituting
a copy of $H$ for every vertex of $G$, is formally defined by taking a disjoint copy of $H$, $H_v$,
for every vertex $v$ of $G$, and joining every vertex in $H_u$ to every vertex in $H_v$ if and only
if $u$ is adjacent to $v$ in $G$.

\begin{teorem} {\rm(\cite{ISRN,brown})} \label{hin}
For any graph $G$, $D(G\diamond K_t, x) = D(G, (1 + x)^t-1)$.
\end{teorem}

\nt We shall show an application of Theorem \ref{hin} to the Turing complexity of the domination
polynomial. If reader like to learn more about the basics of counting complexity theory, can
refer to \cite{complex}. Since the number of dominating sets of graph $G$, i.e. $D(G,1)$ is \#P-complete,
computing domination polynomial with respect to Turing reductions is \#P-hard, even for
restricted graph classes, see e.g. \cite{kij}.

\nt
The value of a graph polynomial at a specific point can give sometimes a surprising information
about the structure of the graph, see e.g. \cite{gcom}. Hardness of computation of graph
polynomial at specific points is another step towards understanding the complexity of a
particular graph polynomial. For example in \cite{Ja}, has shown that the Tutte polynomial
is \#P-hard to compute for any rational evaluation, except those in a semi-algebraic set of
low dimension which are polynomial-time computable.

\nt Let to denote by $D(-,\lambda)$ the problem of computing for an input graph $G$ the evaluation
$D(G,\lambda)$ of the domination polynomial.

\begin{teorem}
 The computation of parameter $D(G,\lambda)$ is \#P-hard, for every $\lambda \in  \mathbb{Q} \setminus \{-2,-1,0\}$.
\end{teorem}

\nt{\bf Proof.} Suppose that $\lambda \in  \mathbb{Q} \setminus \{-2,-1,0\}$. We present an algorithm such that for an input
graph $G$ of order $n$, computes $D(G,x)$ in polynomial time in $n$ using an oracle to $D(-,\lambda)$.
Since $D(G, x)$ is \#P-hard, $D(-,\lambda)$ is \#P-hard. The algorithm is:

\begin{enumerate}
\item[(i)]
 For every $t \in \{1,\ldots, n + 1\}$, compute $D(G\diamond K_t,\lambda) = D(G,(1 +\lambda)^t-1)$.

$D(G\diamond K_t, \lambda)$ is computed using the oracle to $D(-,\lambda)$. Therefore by Theorem \ref{hin},
$D(G, (1 +\lambda)^t-1)$ is computed.

\item[(ii)]
 Interpolate $D(G, x)$ from the values
$(x_0,D(G, x_0)) = ((1 + \lambda)^i- 1,D(G,(1 +\lambda)^i- 1))$,
for $i = 1,\ldots, n + 1$. Since the values $(1 +\lambda)^r-1$ are pairwise distinct (note that $\lambda \not\in  \{-2,-1,0\}$) and $D(G,x)$ has degree $n$, $D(G,x)$ can be interpolated from the computed
values.\quad\qed

\end{enumerate}

\section{Graphs with domination roots in the right half-plane}

\nt The roots of domination polynomial was studied recently by
several authors, see \cite{contem,euro,brown}.

\nt It is clear that $(0,\infty)$ is zero free interval for
domination polynomial of a graph. It is interesting that to
investigate graphs which have complex domination roots with
positive real parts.

\begin{figure}[!h]
\hspace{1.6cm}
\includegraphics[width=8.5cm,height=2.cm]{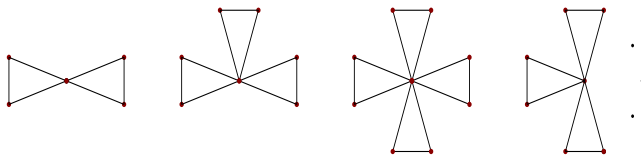}
\caption{\label{Dutch} Dutch-Windmill graphs $G_3^2,G_3^3,G_3^4$ and $G_3^n$, respectively. }
\end{figure}

\nt We consider the graphs obtained by selecting one vertex in each of $n$ triangles and identifying
them. Some call them Dutch Windmill Graphs \cite{site}. See Figure \ref{Dutch}. We denote
these graphs by $G_3^n$. Note that these graphs also called friendship graphs.

\nt We obtain the domination
polynomial of theses graphs and show that there are some of these
graphs whose have complex domination roots with positive real
parts.

\nt We need some
preliminaries.

\begin{teorem}{\rm(\cite{ars})} \label{theorem1}
For every $n\in \mathbb{N}$
$$D(K_{1,n}, x) = x^n + x(1 + x)^n.$$
\end{teorem}

\begin{teorem}{\rm(\cite{brown})}
 The domination polynomial of the star graph, $D(K_{1,n}, x)$, where $n\in \mathbb{N}$, has a real root in the interval
 $(-2n,-ln(n))$, for $n$ sufficiently large.
 \end{teorem}

\begin{figure}[h]
\hspace{1.3cm}
\includegraphics[width=9cm, height=7cm]{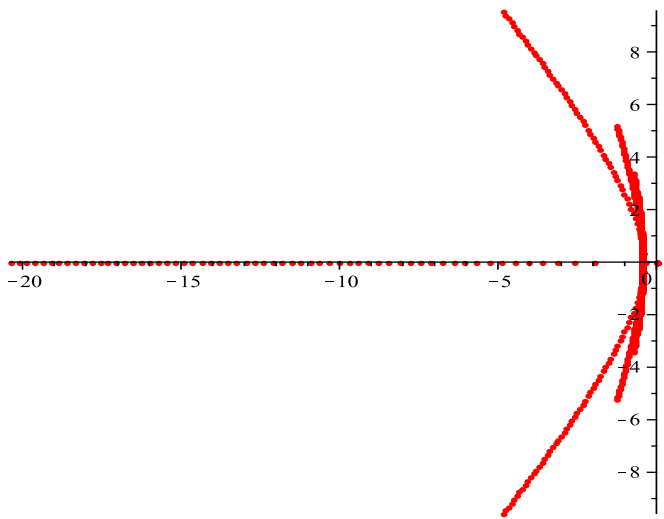}
\caption{\label{figure2} The domination roots of $K_{1,n}$ for $1 \leq n\leq 60$. }
\end{figure}

\nt The domination roots of $K_{1,n}$ for $1 \leq n\leq 60$ has shown in Figure \ref{figure2}.

\nt The join $G = G_1 + G_2$ of two graph $G_1$ and $G_2$ with disjoint vertex sets $V_1$ and $V_2$ and
edge sets $E_1$ and $E_2$ is the graph union $G_1\cup G_2$ together with all the edges joining $V_1$ and
$V_2$.

\begin{teorem}\label{theorem7}{\rm(\cite{euro})}
Let $G_1$ and $G_2$ be  graphs of orders $n_1$ and $n_2$,
respectively. Then
\[
D(G_1+
G_2,x)=\Big((1+x)^{n_1}-1\Big)\Big((1+x)^{n_2}-1\Big)+D(G_1,x)+D(G_2,x).
\]
\end{teorem}

\begin{teorem} \label{d.p.dutch}
For every $n\in \mathbb{N}$,
$$D(G_3^n,x) = (2x + x^2)^n +x(1 + x)^{2n}.$$
\end{teorem}
\nt {\bf Proof.}
It is easy to see that $G_3^n$ is join of $K_1$ and $nK_2$. Now by Theorem \ref{theorem7} we have the result.\quad\qed

\nt In \cite{contem} the following problem has stated:

\nt{\bf Problem.} Characterize all graphs with no real domination
root except zero.

\nt One of the family with no nonzero real domination roots is
$K_{n,n}$ for even $n$:

\begin{teorem}\label{theorem6.2.2}
For every even $n$, no nonzero real numbers is  domination root of
$K_{n,n}$.
\end{teorem}
\nt{\bf Proof.}  It is easy to see that
\[
D(K_{n,n},x)=\Big((1+x)^n-1\Big)^2+2x^n.
\]
If $D(K_{n,n},x)=0$, then $\Big((1+x)^n-1\Big)^2=-2x^n$.
 Obviously this equation does not have real nonzero solution for even $n$.\quad\qed

\nt Domination roots of complete bipartite graphs have been studied extensively in \cite{brown}. We
need the following definition to state one of the main result on domination roots of $K_{n,n}$.

\begin{define}  If $\{f_n(x)\}$ is a family of (complex) polynomials, we say that a number
$z\in  \mathbb{C}$ is a limit of roots of $\{f_n(x)\}$ if either $f_n(z) = 0$ for all sufficiently large $n$ or $z$ is a
limit point of the set $R(f_n(x))$, where $R(f_n(x))$ is the union of the roots of the $f_n(x)$.
\end{define}

\begin{figure}[h]
\hspace{2.3cm}
\includegraphics[width=6cm, height=6.5cm]{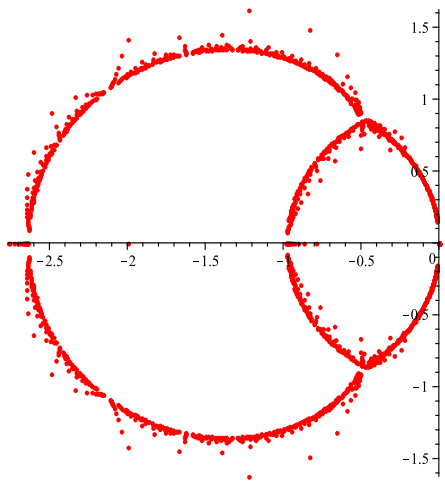}
\caption{\label{figure3} Domination roots of $K_{n,n}$ for $1 \leq n\leq 40$. }
\end{figure}

\nt The domination roots of $K_{n,n}$ for $1\leq n \leq 40$ has shown in Figure \ref{figure3}. See also \cite{brown}. As we
can see the domination roots of $K_{n,n}$ are bounded. The following theorem characterize
limit of roots of the domination polynomials of $K_{n,n}$ for every $n\in  \mathbb{N}$.

\begin{teorem}{\rm(\cite{brown})}
The complex numbers $z$ that satisfy any of the following conditions:
\begin{enumerate}
\item[(i)] $|z-(-1)| = 1$, $\mathfrak{R}(z) > \frac{-1}{2}$,

\item[(ii)] $z = \frac{-1}{2}\pm \frac{\sqrt{3}}{2}i$,

\item[(iii)] $ |1 + z|^2 = |z|$, $\mathfrak{R}(z) < \frac{-1}{2}$,

\end{enumerate}
are limits of roots of the domination polynomial of the graphs $K_{n,n}$,  $n\in  \mathbb{N}$.

\end{teorem}

\nt Here we prove that, for every odd natural number $n$,  Dutch windmill graph $G_3^{n}$ have no real
roots except zero.

\begin{teorem}
For every odd natural number $n$, no nonzero real numbers is  domination
root of $G_{n}^n$.
\end{teorem}
\nt{\bf Proof.} By Theorem \ref{d.p.dutch}, for every $n\in
\mathbb{N}$, $D(G_3^n,x) = (2x + x^2)^n +x(1 + x)^{2n}.$ If
$D(G_3^n,x)=0$, then we have
\[
x=-\big(1-\frac{1}{(1+x)^2}\big)^n.
\]
First suppose that $x\geq 0$. Obviously the above equality is true
just for real number 0, since for nonzero real number the left
side of equality is positive but the right side is negative. Now
suppose that $x < -2$. In this case the left side is less than $-2$
and the right side $-\big(1-\frac{1}{(1+x)^2}\big)^n$ is greater
than $-1$, a contradiction. Finally we shall consider $-2<x<0$.
 In this case obviously the above equality is not true for any real number. Because for odd $n$ and  for every real number $-2<x<0$, the left side of equality is negative but the right side is positive. \quad\qed

\begin{figure}[h]
\begin{minipage}{5.9cm}
\includegraphics[width=\textwidth]{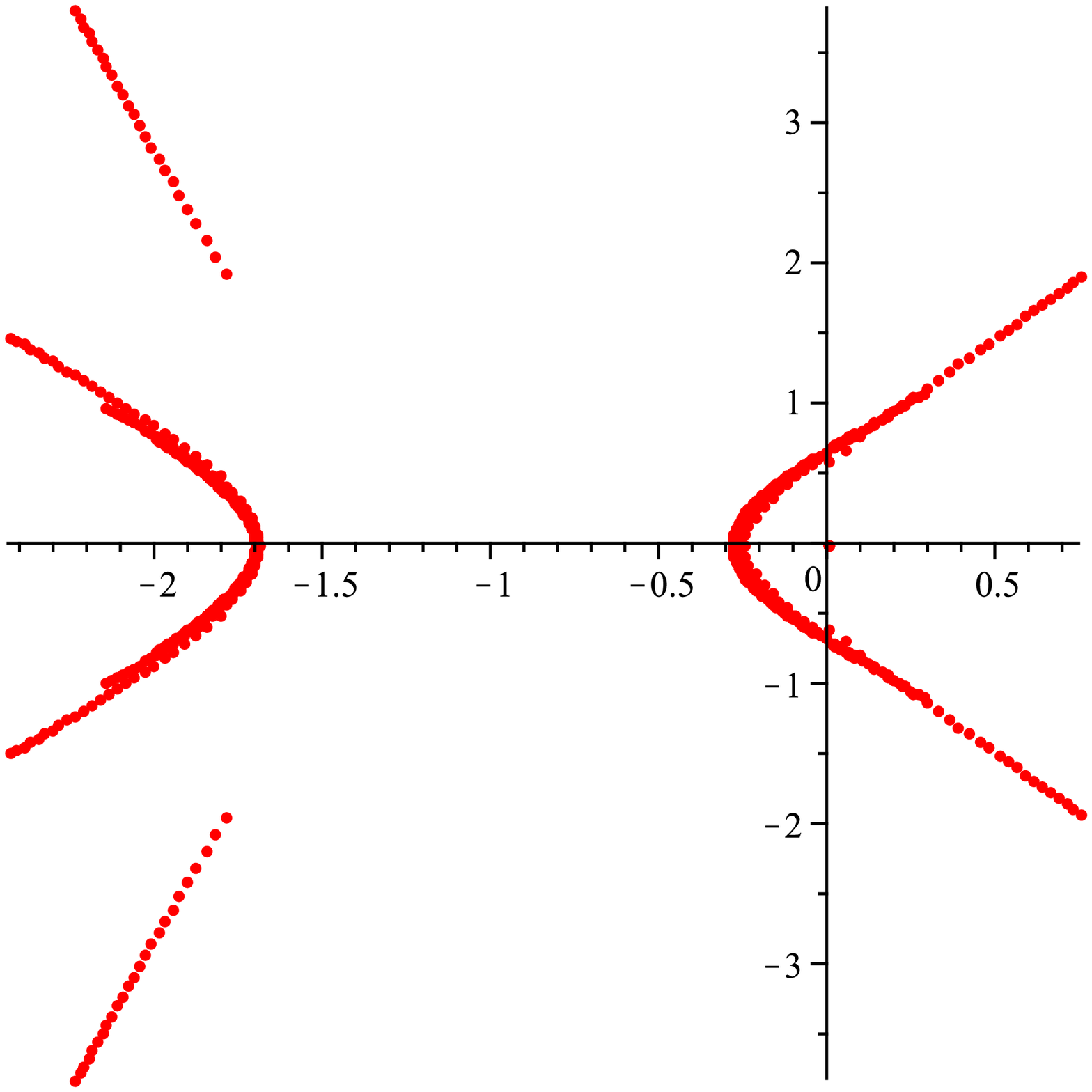}
\end{minipage}
\begin{minipage}{5.9cm}
\includegraphics[width=\textwidth]{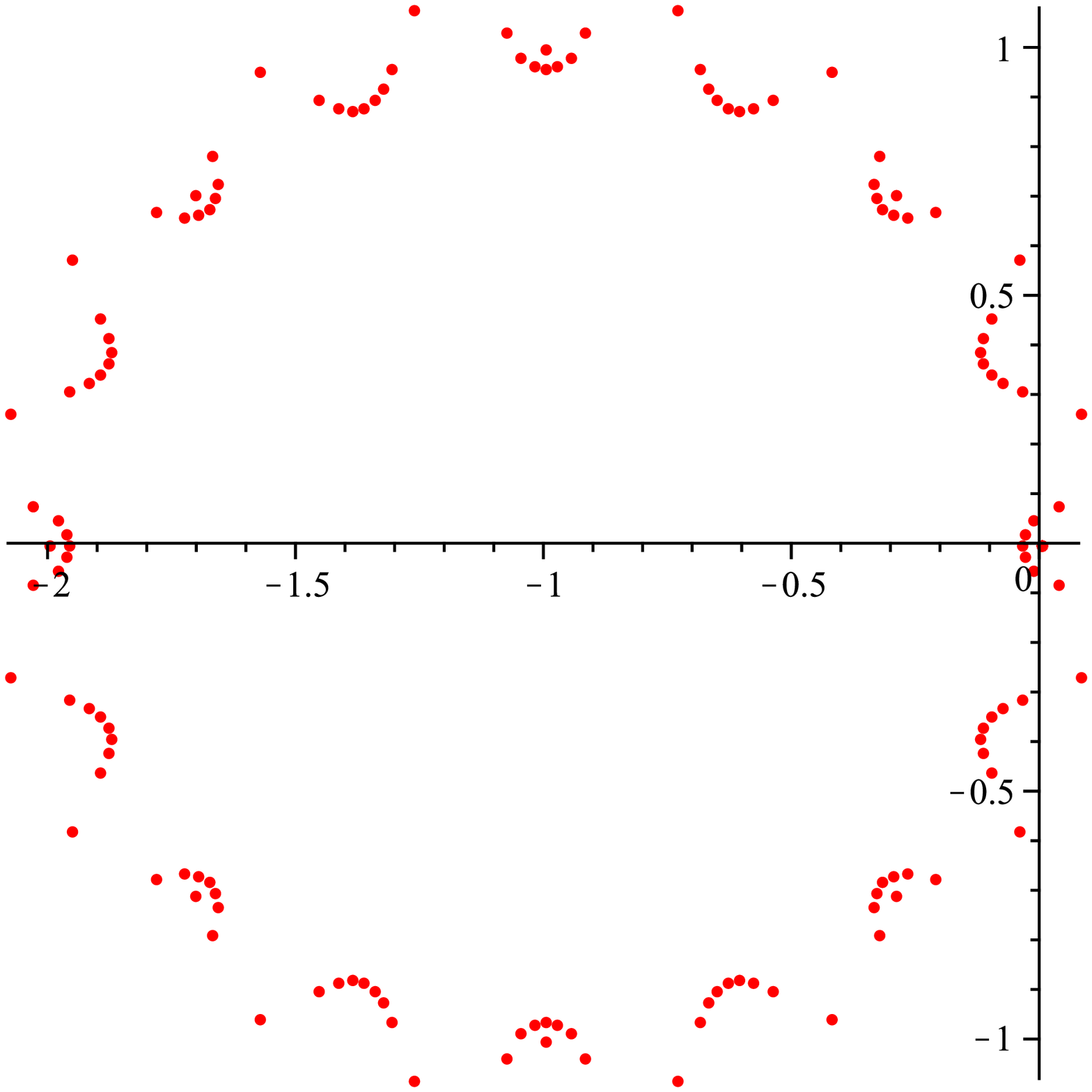}
\end{minipage}
\caption{\label{figure2'} Domination roots of graphs $G_3^n$ for $1 \leq n \leq 30$, and $G_3^{8}\diamond K_{8}$,  respectively.}
\end{figure}

\nt{\bf Remark.} Using Maple we observed that the domination
polynomial of $G_3^n$ for $n\geq 6$  have complex  roots with
positive real  parts. For example $D(G_3^6,x)$ has complex root
with real part $0.0003550296365$. See Figure \ref{figure2'}.

\nt Using Theorems \ref{hin} and \ref{d.p.dutch} we have the following theorem:

\begin{teorem}
$D(G_3^n \diamond K_t, x) =((1 + x)^{2t}- 1)^n +
((1 + x)^t- 1)(1 + x)^{2nt}.$
\end{teorem}

\nt It is interesting that the families of graphs $G_3^n\diamond K_t$ have domination roots with positive
real parts (see Figure \ref{figure2'}).


\section{The domination polynomial of other classes of graphs}

\nt In this section we study the domination polynomial of other classes of graphs.

\nt The vertex contraction $G/v$ of a graph $G$ by a vertex $v$ is
 the operation under which all vertices in $N(v)$ are joined to
each other and then $v$ is deleted (see \cite{walsh}).

\begin{teorem}  \label{theorem2}{\rm(\cite{book,TMU,Kot})}
For any vertex $v$ in a graph $G$ we have $$D(G,x) = xD(G/v,x) +
D(G-v,x) + xD(G-N[v],x)-(x + 1)p_v(G,x)$$ where $p_v(G,x)$ is the
polynomial counting those dominating sets for $G-N[v]$ which
additionally dominate the vertices of $N(v)$ in $G$.
\end{teorem}

\nt Theorem \ref{theorem2} can be used to give a recurrence
relation which removes triangles. We define a new operation on
edges incident to a vertex $u$: we denote by $G\odot u$ the graph
obtained from $G$ by the removal of all edges between any pair of
neighbors of $u$. Note $u$ is not removed from the graph. The
following  recurrence relation is useful on graphs which have many
triangles. This following result also appear in \cite{Kot} but were proved independently.

\begin{teorem}\label{theorem3} 
 Let $G = (V,E)$ be a graph and $u \in V$. Then
\[
D(G, x) = D(G - u, x) + D(G \odot u, x) - D(G\odot u - u, x)
\]
\end{teorem}

\nt{\bf Proof.} Since  the operation $\odot u$ only removes the
edges between vertices in $N(u)$, we have the following relations:

$(G \odot u)/u \cong G/u$, $p_u(G,x) = p_u(G\odot u, x)$, $(G
\odot u)- N[u]\cong (G- N[u]).$

\nt Using these relations and Theorem \ref{theorem2}, we have
\begin{eqnarray*}
D(G, x)- D(G-u, x) &=& xD(G/u, x) + xD(G- N[u], x) -(1+x)p_u(G, x)\\
& =& xD((G \odot u)/u,x) + xD((G\odot u)- N[u],x)\\
& -&(1 +x)p_u(G\odot u,x)
\end{eqnarray*}

 \nt Now by   Theorem \ref{theorem2} for  $G\odot u$ we have
 \begin{eqnarray*}
xD((G \odot u)/u,x) &+& xD((G\odot u)- N[u],x)-(1 +x)p_u(G\odot u,x)\\
&=&D(G \odot u,x) - D((G\odot u)-u,x)
\end{eqnarray*}
 Therefore we have the result.\quad\qed

\nt Now we use Theorem \ref{theorem3} to study the domination polynomial and domination roots of
some other classes of graphs:

\begin{enumerate}

\item[(1)]  The Dutch Windmill graph with an extra edge $vu$, i.e pendant edge to central vertex. The
3 graphs $G-u$, $G\odot u$ and $G\odot u-u$ are (i) $K_2$,   $n$
times and $K_1$; (ii) $K_{1,2n+1}$;  (iii) $P_1$,  $2n+1$
times, respectively. So, by Theorem \ref{theorem3} we have
\begin{eqnarray*}
           D(G_n, x) &=& x(2x+x^2)^n  + x^{2n+1} + x(1+x)^{2n+1} - x^{2n+1}\\
           & =&x\big((x^2+2x)^n+(x+1)^{2n+1}\big).
           \end{eqnarray*}

\nt The reader is able to see the sequence of coefficients of this polynomial in the site
``The on-line encyclopedia of integer sequences" (\cite{seq}) as A213658.

\nt It is interesting that the graph $G_{n}$ has domination roots in the right-half plane
(Figure \ref{figure5}). Also for these kind of family of graphs we have the following theorem:

\begin{teorem}
 For every natural $n$, there is exactly one nonzero real domination
root of $G_n$.
\end{teorem}
\nt {\bf Proof.} Suppose that $\alpha$ is a root of $D(G_n,x)$. So we have
$$(1 + \frac{1}{(\alpha+1)^2-1})^n = \frac{-1}{\alpha+1}.$$
By substituting $\alpha + 1 = t$, we shall have $t^{2n+1} + (t^2 - 1)^n = 0$. This equation has
only one real root in $(0, 1)$ for odd $n$, and has only one real root in $(-1, 0)$ for even
$n$. Therefore $D(G_n,x)$ has only one real root in $(-1,0)$ or in $(-2,-1)$.\quad\qed

\nt As you can see in the Figure \ref{figure5} the graph $G_{n}$ has domination roots in the right-half
plane.

\item[(2)] The fan graph. A fan graph $F_{m,n}$ is defined as the
graph join $\overline{K_m}+P_n$, where $\overline{K_m}$ is the
empty graph on $m$ vertices  and is $P_n$ the path graph on $n$
vertices.
Here we  consider $F_{2,n}$.
To obtain $F_{2,n}$ take two vertices $u,~v$ and join each of $n$
vertices $1,2,...,n$ to both $u$ and $v$.
By Theorem \ref{theorem7}, we have
\begin{eqnarray*}
D(F_{2,n},x) &=& (x^2+2x)((1+x)^n-1) + x^2 +D(P_n,x) \\
& =&2x((1+x)^n-1) +D(P_n,x).
\end{eqnarray*}

\nt We can see the sequence of coefficients of this polynomial in the site ``The on-line
encyclopedia of integer sequences" (\cite{seq}) as A213657.

\nt The domination roots of graph $F_{2,n}$ has shown in Figure \ref{figure5}.

\begin{figure}[h]
\begin{minipage}{5.9cm}
\includegraphics[width=\textwidth]{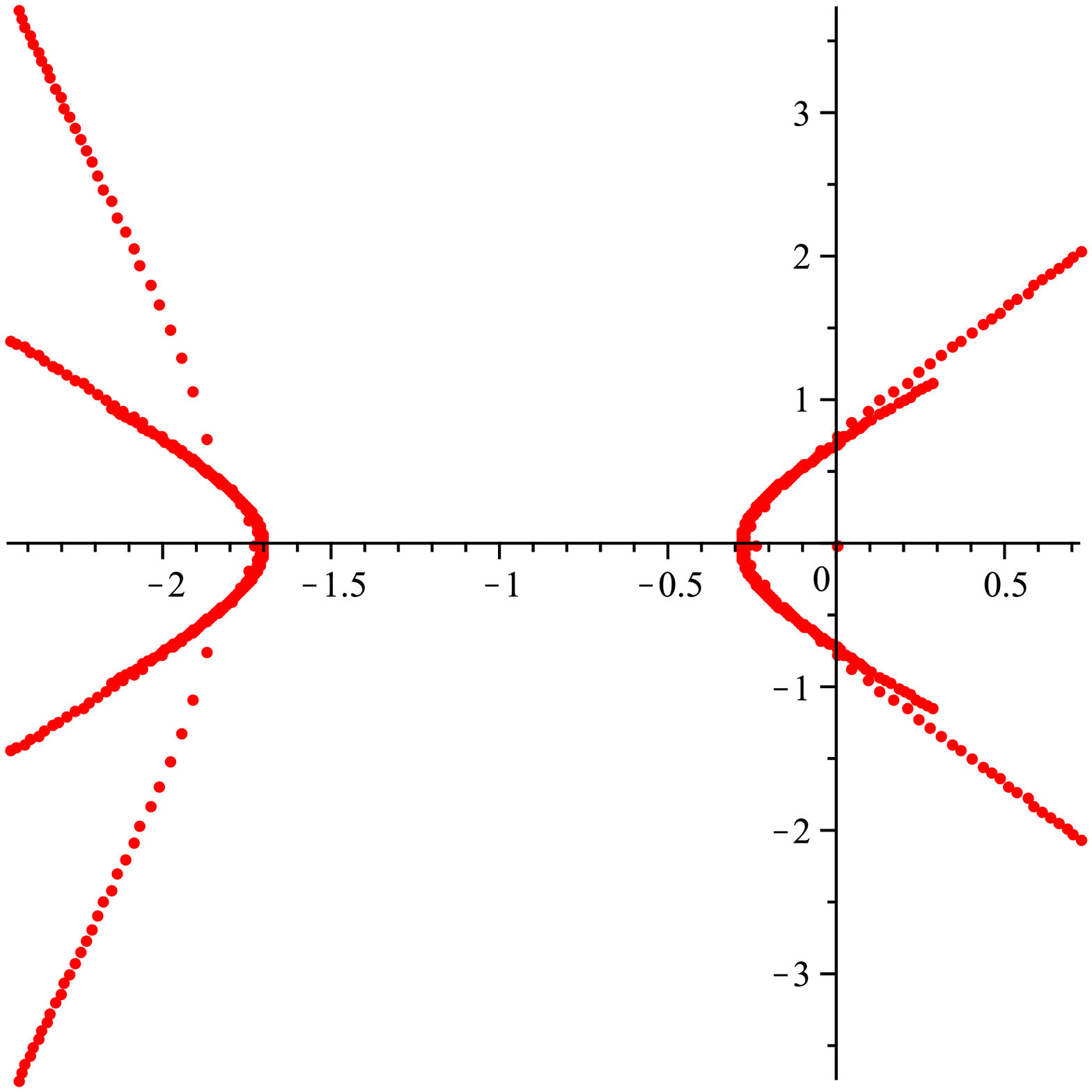}
\end{minipage}
\begin{minipage}{5.9cm}
\includegraphics[width=\textwidth]{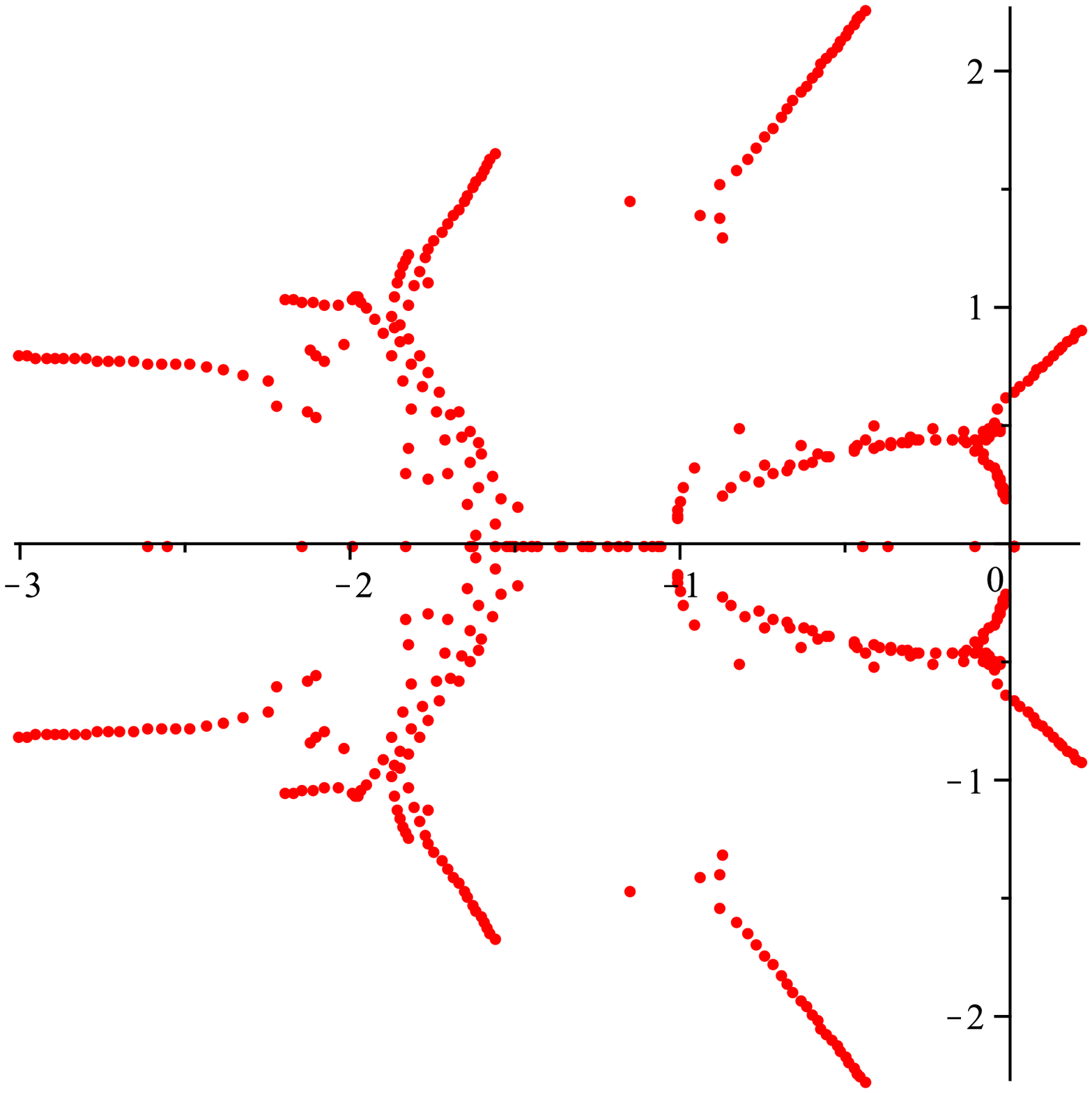}
\end{minipage}
\caption{\label{figure5} Domination roots of graphs $G_{n}$ and $F_{2,n}$ for $1 \leq n\leq 30$, respectively.}
\end{figure}


\nt
 Using Theorem \ref{theorem7}, we have the following corollary:

 \begin{korolari}
For every natural $m,n \in \mathbb{N}$,
 $$D(F_{m,n},x) = ((1+x)^m-1)((1+x)^n-1)+x^m+D(P_n,x).$$
           \end{korolari}

\item[(3)] The Gem graph $G$. Consider the path $P_{n+1}$ and an
additional vertex $u$; join $u$ to each vertex of the path. The 3
graphs in Theorem \ref{theorem3}  are: (i) $P_{n+1}$; (ii) the
star $K_{1,n+1}$; (iii) $P_1$,  $(n+1)$ times. So,
\begin{eqnarray*}
           D(G,x) &=& D(P_{n+1},x) + x^{n+1} + x(1+x)^{n+1} - x^{n+1} \\
           &=&D(P_{n+1},x)+ x(1+x)^{n+1}.
           \end{eqnarray*}
\nt The reader is able to see the sequence of coefficients of this polynomial in the site
``The on-line encyclopedia of integer sequences" (\cite{seq}) as A213662.

\item[(4)]  The Gem graph  with an extra edge $vu$ which is
denoted by $G'$. The 3 graphs in Theorem \ref{theorem3} are (i)
$P_{n+1}$ and $P_1$; (ii) $K_{1,n+2}$; (iii) $P_1$,  $(n+2)$
times. So,
\begin{eqnarray*}
           D(G',x)& =& xD(P_{n+1},x) + x^{n+2} + x(1+x)^{n+2} - x^{n+2} \\
           &=& x\big(D(P_{n+1},x)+(1+x)^{n+2}\big).
           \end{eqnarray*}

\item[(5)]  Join a vertex $u$ with two consecutive vertices of the
cycle $C_n$ (i.e. a triangle placed on an edge of $C_n$). Let to
denote this graph by $G$. The 3 graphs in Theorem \ref{theorem3}
are: (i) $C_n$; (ii) $C_{n+1}$; (iii) $P_n$. So,
          $$D(G,x) = D(C_n,x) + D(C_{n+1},x) - D(P_n, x).$$

\nt We can see the sequence of coefficients of this polynomial in the site ``The on-line
encyclopedia of integer sequences" (\cite{seq}) as A213664.

\item[(6)]  The wheel graph $W_n$. The 3 graphs in Theorem
\ref{theorem3} are (i) $C(n-1)$; (ii) $K_{1,n-1}$; (iii) $P_1$,
$n-1$ times. So,
\begin{eqnarray*}
         D(W_n,x) &=& D(C_{n-1},x) + x^{n-1} + x(1+x)^{n-1} - x^{n-1} \\
         &=&D(C_{n-1},x)+x(1+x)^{n-1}.
         \end{eqnarray*}

         \end{enumerate}

\begin{figure}[h]
\begin{minipage}{5.9cm}
\includegraphics[width=\textwidth]{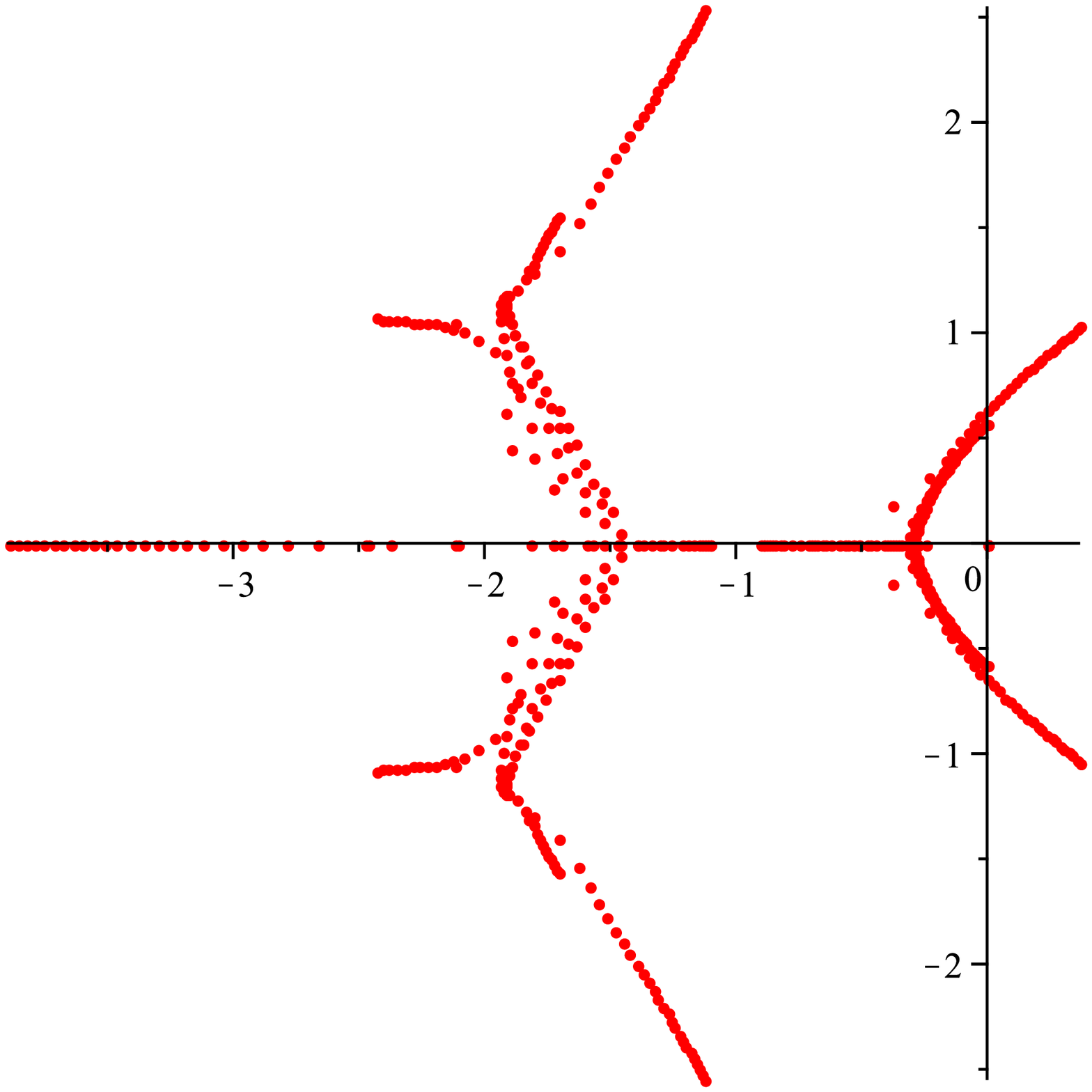}
\end{minipage}
\begin{minipage}{5.9cm}
\includegraphics[width=\textwidth]{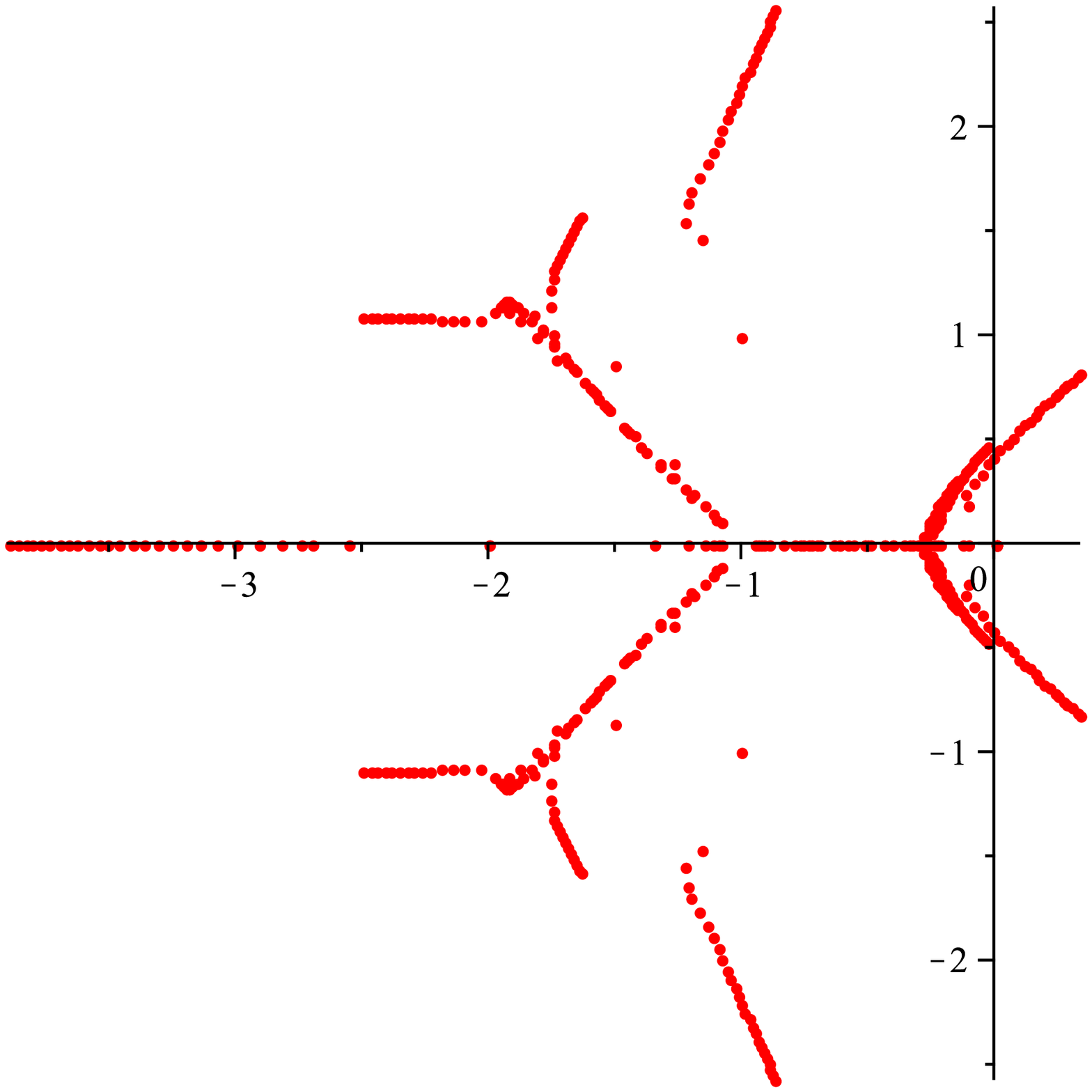}
\end{minipage}
\caption{\label{figure6} Domination roots of graphs $G'_{n}$ and $W_{n}$ for $1 \leq n\leq 30$, respectively.}
\end{figure}

\nt As we can see in the Figure \ref{figure6}, there are graphs in the families of $G'_n$ and $W_n$ which their
domination roots are in the right half-plane.

\nt {\bf Acknowledgement.} The authors would like to express their  gratitude to the referee for helpful comments. The research of the first  author was in part supported by a grant from IPM (No. 91050015)
and partially supported by Yazd University Research Council. The authors would like to thank S. Jahari for some useful Maple
programme.



\end{document}